\documentclass[draft]{article}
\usepackage[T2A]{fontenc}
\usepackage[cp1251]{inputenc}
\usepackage{amssymb,amsmath,amsthm,amsfonts}
\begin{document}
%%\addcontentsline{toc}{section}{Хабибуллин~Б.~Н., Таминдарова~Н.~А. Распределение нулей голоморфных функций: условия типа Бляшке}
%Для содержания на русском языке
\begin{center}
\textbf{ON DISTRIBUTION OF ZEROS\\OF HOLOMORPHIC FUNCTIONS:\\
BLASCHKE-TYPE CONDITIONS}%
\footnote[1]{Supported by RFBR (project 13-01-00030)}\\
{\sc B.~N.~Khabibullin, N.~R.~Tamindarova} (Ufa)\\
{\it Khabib-Bulat@mail.ru}
\end{center}
\textbf{0. Introduction.} We use the definitions and information from [1]. Our results are closely related to [2]--[22].
As usual, $\mathbb N:=\{1,2, \dots\}$, $\mathbb R$ and $\mathbb C$ are the sets of real and complex numbers, resp. 
The symbol $\mathbb C_{\infty} := \mathbb C \cup  \{\infty\}$  denotes the Riemann sphere (the extended complex plane);  $D (z,r):=\bigl\{ z' \in \mathbb C\colon  |z'-z| < r\bigr\}$ for $z\in \mathbb C$ and $r>0$; $\mathbb D := D(0,1)$.
Let $D$ be a subdomain of\/ $\mathbb C_{\infty}\neq D$ always here. We denote by $\text{Hol}(D)$ and $\text{sbh}(D)$ the class of holomorphic  and subharmonic functions on $D$, resp. The class $\text{sbh}(D)$  contains the function $\boldsymbol{-\infty}\colon z\mapsto -\infty$, $z\in D$. 
An elementary consequence of classical Nevalinna Theorem is

\smallskip
\noindent
	\textsc{Theorem} (with W.~Blaschke condition). {\it  Let  $M\in \text{\rm sbh}(\mathbb D)\setminus \{\boldsymbol{-\infty}\}$ with the  Riesz measure $\nu_M$ 
	and a version of the Blaschke condition  
				\begin{equation*}\label{Bcnu}
\int\limits_{1/2\leq |z|<1} \log \frac1{|z|}\, d \nu_M(z)<+\infty
		\eqno{\rm (Bm)}
		\end{equation*}
	is fulfilled. If  a non-zero function $f\in \text{\rm Hol} (\mathbb D)$ such that  $\log |f|\leq M$ on $\mathbb D\setminus D(0,1/2)$ (pointwise) and $f$ vanish on a sequence\/ ${\tt Z}=\{{\tt z}_k\}_{k=1,2,\dots}\subset \mathbb D$ (we wrie $f({\tt Z})=0$), then the Blaschke condition   
		\begin{equation*}\label{BcZ}
		\sum_{1/2\leq |{\tt z}_k|<1} \log \frac1{|{\tt z}_k|}<+\infty	
		\eqno{\rm (Bz)}
	\end{equation*}
also  is fulfilled.} 
Here the test function $z\mapsto \log \dfrac{1}{|z|}$, $z\in \mathbb D$, from (Bm)--(Bz) is positive (sub)har\-m\-o\-n\-ic function on $\mathbb D\setminus D(0,1/2) $, 
and 
	$$
	\lim\limits_{1>|z|\to 1} \log \frac1{|z|} =0.
	\eqno{(0)}
	$$
 For a subset $S$ of $\mathbb C_{\infty}$   we denote by $\overline{S}$ and $\partial S$ the closure of $S$ and the boundary of $S$ in $\mathbb C_{\infty}$.
 If $\overline{S}$ is a compact subset of $D$ in the topology induced from $\mathbb C_{\infty}$, then  we write $S \Subset D$. For $S\subset \mathbb C_{\infty}$ we denote by $\text{sbh} (S)$ the class  of functions that  are subharmonic on some open set containing $S$; $\text{sbh}^+(S):=\bigr\{u\in \text{sbh}(S)\colon u\geq 0 \text{  on $S$}\bigl\}$.

\smallskip
\noindent
{\sc Main Problem.} {\it Let $D_0\Subset D$ be a subdomain, $M\in \text{\rm sbh} (D)$ with the Riesz measure $\nu_M$, $ f\in \text{\rm Hol} (D)$, 
  ${\tt Z}=\{{\tt z}_k\}_{k=1,2,\dots}\subset D$, $f({\tt Z})=0$, and $\log |f|\leq M$ on $D\setminus D_0$. 
	Under what conditions on the   function $v\colon D\setminus D_0\to [0,+\infty)$, for  $f\not\equiv 0$,  the implication
	\begin{equation*}
	\left( \, \int\limits_{D\setminus D_0}v\, d \nu_M <+\infty\right) \; \Longrightarrow \; \left(\, \sum_{{\tt z}_k \in D\setminus D_0}v ({\tt z}_k) <+\infty\right)
\eqno{\text{$(\Rightarrow)$}}
\end{equation*}
is true}\,? This gives uniqueness theorems: if a series from  $(\Rightarrow)$ diverges, the integral from $(\Rightarrow)$ is finite, 
and $\log |f|\leq M$ on $D\setminus D_0$, then  $f\equiv 0$.

\smallskip
\noindent
{\sc Our Solution} (see Corollary 1 and Theorem 1 below). {\it  It is sufficient that  $v\in \text{\rm sbh}^+(D\setminus D_0)$, and\/ {\rm (cf. (0))}}
\begin{equation*}
	\lim\limits_{z\to \partial D} v=0,
	\eqno{\rm (O)}
\end{equation*}
i.\,e. for each $\varepsilon >0$ there is  $ D_{\varepsilon}\Subset D$ such that  $|v|<\varepsilon$ on $D\setminus (D_{\varepsilon}\cup D_0)$.
\smallskip
  
\noindent
1. \textbf{$\delta$-subharmonic functions} [23]--[24].  For a function or a mesure $a$ we denote by $a\bigm|_S$ the restriction of $a$ to $S$.   $L^1_{\text{loc}} (D)$ denotes  the set of locally integrable functions $f$ on $D$ with respect to the restriction $\lambda\bigm|_D$ of Lebesgue measure $\lambda$. We 
also consider the function  $\boldsymbol{+\infty}\colon z\mapsto +\infty$, $z\in D$.  Two  functions $\boldsymbol{\pm\infty}$ are $\delta$-subharmonic functions.  
 Another function  $M\colon D\to [-\infty,+\infty]$ is called $\delta$-subharmonic on $D$ if $M\in L^1_{\text{loc}} (D)$ and
\begin{enumerate}
\item For any subdomain  $D_0\Subset D$ there exists a constant $C_0\in [0,+\infty )$ such that for each finite  
 function $\phi \colon D\to \mathbb R$ with $\text{supp\,} \phi \subset D_0$, $\phi \in C^2(D)$, the inequality
			\begin{equation*}
			\left|\, \int\limits_{D_0} M	\Delta\phi \, d\lambda \right|\leq C_0 \max_{z\in D_0} \bigl|\phi (z)\bigr|
			\end{equation*}
			is fulfilled.  Further  $\nu_M:=\frac{1}{2\pi} \Delta M$ is the Riesz charge of $M$ where $\Delta$ is the Laplace operator acting in the sense of distribution theory. Besides,  we have the Hahn--Jordan decomposition $\nu_M:=\nu_M^+-\nu_M^-$ where  $\nu_M^+$ and $\nu_M^-$ are called the positive and negative part of $\nu_M$; $|\nu_M|:=\nu_M^++ \nu_M^-$ is the absolute variation of $\nu_M$.
	\item		We define $\text{dom}_M\subset D$ as the set of points $z\in D$ such that there is $r_z>0$ 
	  with the properties 	$D(z,r_z)\Subset D$ and
	\begin{equation*}
	\int\limits_{D(z,r_z)} \log |z'-z_0| \, d |\nu_M|(z')>-\infty.
		\end{equation*}
		For $z\in \text{dom}_M$ we set 
		\begin{equation*}
			M(z)=\lim_{0<r\to 0} \frac{1}{\pi r^2} \int\limits_{D(z,r)} M \, d \lambda.
			\end{equation*}
			\item Here it is convenient to set $M(z)=+\infty$ for all $z\in D\setminus \text{dom}_M$.
\end{enumerate}
We denote by $\delta\text{-sbh}(D)$ the class of $\delta$-subharmonic functions on $D$. 
Each  function $M\in \delta\text{-sbh}(D)\setminus\{\boldsymbol{\pm\infty}\}$   admits a unique representation
$M=M_+-M_-$ where $M_+, M_-\in \text{sbh\,}(D)$ with the Riesz measures $\nu_M^+,\nu_M^-$, resp.

\noindent
 \textbf{2. Main results.}
For a regular domain $D \subset  \mathbb C_{\infty}$ we denote by   $g_D (\cdot,z_0)$ the extended Green’s function of $D$ with pole at $z_0\in  D$, that is, $g_D (z',z_0) \equiv 0$ for all $z'\in \mathbb C \setminus  D$, $g(\cdot,z_0)\bigm|_{D\setminus\{z_0\}}$ is harmonic, and $g(\cdot,z_0)$ is subharmonic on  $\mathbb C_{\infty} \setminus \{z_0\}$.

We denote by  $\text{const}_{a_1, a_2, \dots}$  a constant depending only on $a_1,a_2, \dots$.

Let\/ $D_0, D$ are domains in\/ $\mathbb C_{\infty}$,  and  $\varnothing\neq  D_0\Subset  D\neq \mathbb C_{\infty}$.

Let $b\in [0,+\infty)$. We set $ \text{sbh}_0^+(D\setminus D_0;\leq b)$
\begin{equation*}
	:=\left\{v\in \text{sbh}^+(D\setminus D_0)\colon 
	\lim_{z\to \partial D} v(z)\overset{\text{see (O)}}{=}0, \; \sup\limits_{z\in \partial D_0}v(z)\leq b \right\}.
\end{equation*}
We consider a function  $M\in \delta\text{-sbh\,}(D)\setminus \{\boldsymbol{\pm\infty}\}$ with the Riesz charge $\nu_M$.

\smallskip
\noindent
\textsc{Main Theorem.} {\it FOR ANY 
\begin{enumerate}
	\item[\rm (i)] point $z_0\in D_0\cap \text{\rm dom}_M$ and number $b\in [0,+\infty)$,
			\item[\rm(ii)] regular domain $\widetilde{D}\subset \mathbb C_{\infty}$, $D_0\Subset \widetilde{D}\subset D$, $\mathbb C_{\infty}\setminus 
		\overline{\widetilde{D}}\neq \varnothing$, 
	\end{enumerate}
	THERE EXISTS a number
$C:=\text{\rm const}_{z_0,D_0,\widetilde{D},D,b}>0$,
such that for each $u\in \text{\rm sbh\,} (D)\setminus \{\boldsymbol{-\infty}\}$ the inequality $u\leq M$ on $D\setminus D_0$ entails that,
for any function $v{\in}  \text{\rm sbh}_0^+(D\setminus D_0;\leq b) $, the inequality
\begin{multline*}
C u(z_0) 	+\int\limits_{D\setminus D_0}  v \,d {\nu}_u 		\leq	\int\limits_{D\setminus D_0}  v \,d {\nu}_M	+\int\limits_{\widetilde{D}\setminus D_0}
v \,d {\nu}_M^- \\	+		C\int\limits_{\widetilde{D}} g_{\widetilde{D}}(\cdot, z_0)  \,d {\nu}_M  
		+C\int\limits_{\widetilde{D}\setminus D_0} g_{\widetilde{D}}(\cdot, z_0)  \,d {\nu}_M^-  +CM(z_0)
\end{multline*}
is fulfilled. Moreover, if $\widetilde{D}\Subset D$, then this estimate can be replaced by
\begin{equation*}
	\int\limits_{D\setminus D_0}  v \,d {\nu}_u 		\leq	\int\limits_{D\setminus D_0}  v \,d {\nu}_M	+(b+C)\, \overline{C}_M-C u(z_0),
\eqno{\rm (C)}
\end{equation*}
where  a constant $\overline{C}_M:=\text{\rm const}_{z_0,D_0, \widetilde{D}, D, M}<+\infty$ is positively homog\-e\-n\-e\-o\-us of $M$, i.\,e.
 $\overline{C}_{aM}=a\overline{C}_M$ for  $a\in [0,+\infty)$, and   upper semi-additive of $M$, i.\,e.
		$\overline{C}_{M_1+M_2}\leq \overline{C}_{M_1}+\overline{C}_{M_2}$ for suitable $M_1, M_2$.

Conversely, IF the function $M$ is continuous on \underline{regular} domain $D$ and,   for a Borel 
measure $\nu\geq 0$ on $D$ and for  a number $b>0$, there is a constant $C'$ such that	 for each function 
$v{\in}  \text{\rm sbh}_0^+(D\setminus D_0;\leq b) $ the inequality 
\begin{equation*}
	\int\limits_{D\setminus D_0}  v \,d {\nu} 		\leq	\int\limits_{D\setminus D_0}  v \,d {\nu}_M	+C'
\eqno{\rm(C')}
\end{equation*}
is fulfilled,  THEN there exists a function $u\in \text{\rm sbh}(D)\setminus \{\boldsymbol{-\infty}\}$ with the Riesz measure $\nu_u\geq \nu$ such that $u\leq M$
on  $D$, and, for any function $u_0\in \text{\rm sbh\,} (D)$ with the Riesz measure $\nu_{u_0}=\nu$ and for any number $\varepsilon >0$, there exists a non-zero function $f\in \text{\rm Hol\,}(D)$ such that the inequality  
\begin{equation*}
	u_0(z)+\log \bigl|f(z)\bigr|	\leq 	\frac{1}{2\pi}\int_0^{2\pi} M(z+re^{i\theta}) \, d\theta+ \log \frac{\bigl(1+|z|\bigr)}{r}^{1+\varepsilon}
\eqno{\rm (L)}
\end{equation*}
is fulfilled for all $z\in D$ and for all\/ 
\begin{equation*}
	0<r<\min \Bigl\{1+|z|, \text{\rm dist\,}(z,\partial D):=\inf \bigl\{|z'-z|\colon z'\in \partial D\bigr\}\Bigr\}.
\eqno{\rm(d) }
\end{equation*}
}
\smallskip

\noindent
\textbf{Remark.} In addition, if $D=\mathbb C$ or $D$ is a domain  with smooth boundary of class $C^1$ and, for $v\bigm|_{\partial D}:\equiv 0$, 
		there exists the derivative $\frac{\partial v}{\partial \vec{n}_{\text{\rm in}}}\bigm|_{\partial D}\equiv 0$ along the inner normal $\vec{n}_{\text{\rm in}}$ of $\partial D$, then the integral in the right parts of  the inequalities\/ {\rm (C)} and $(C')$  can be replaced by the  integral 		$\int\limits_{D\setminus D_0} M\, d \nu_v$, where $\nu_v$ is the Riesz measure of $v$.

\smallskip
\noindent
\textsc{Theorem 1.} {\it Let  $M\in \text{\rm sbh\,}(D)$,  $b\in [0,+\infty)$. Then 
	there are numbers $C:=\text{\rm const}_{D_0,D,b}>0$,  $\overline{C}_M:=\text{\rm const}_{D_0,  D, M}\geq 0$
such that, for any  $v{\in}  \text{\rm sbh}_0^+(D\setminus D_0;\leq b) $ and for each $u\in \text{\rm sbh\,} (D)\setminus \{\boldsymbol{-\infty}\}$, the inequality $u\leq M$ on $D\setminus D_0$ entails the inequality\/ 
{\rm (C)}. 

Conversely, if,  for a Borel measure $\nu\geq 0$ on \underline{regular} domain  $D$,  there is a constant $C'$ such that	 for each function 
$v{\in}  \text{\rm sbh}_0^+(D\setminus D_0;\leq b) $ the inequality\/ $(C')$ is fulfilled, then 
there exists $u\in \text{\rm sbh}(D)\setminus \{\boldsymbol{-\infty}\}$ with the  Riesz measure $\nu_u\geq \nu$ such that $u\leq M$
on  $D$, and, for any function $u_0\in \text{\rm sbh\,} (D)$ with the Riesz measure $\nu_{u_0}=\nu$ and for any number $\varepsilon >0$, there exists a non-zero function $f\in \text{\rm Hol\,}(D)$  satisfying\/ $\rm (L)$--$\rm (d)$.

The Remark   remains in force.
}
\smallskip

\noindent
\textsc{Corollary 1.} 
{\it Let $D_0\Subset D$ be a subdomain of $D$, $M\in \text{\rm sbh} (D)$ with the Riesz measure $\nu_M$, $ f\in \text{\rm Hol} (D)$, 
  ${\tt Z}=\{{\tt z}_k\}_{k=1,2,\dots}\subset D$, $f({\tt Z})=0$, and $\log |f|\leq M$ on $D\setminus D_0$.
	If $v\in \text{\rm sbh}^+(D\setminus D_0)$, and $	\lim\limits_{z\to \partial D} v\overset{\rm see (O)}{=}0$, then
	 the implication $(\Rightarrow)$ from Main Problem is true.}
\smallskip

\noindent
\textsc{Corollary 2.} {\it Let $D=\mathbb C$ or $D$ be a domain  with smooth boundary of class $C^1$, and $D_0\Subset D$ be a subdomain of $D$.
 Suppose that a continuous function $v\colon \overline{D}\setminus D_0\to \mathbb R$  satisfies the following conditions
\begin{itemize}
	\item $v\bigm|_{D\setminus D_0}\in \text{\rm sbh}^+(D\setminus D_0)$ with the Riesz mesure $\nu_v$, and $	 v\bigm|_{\partial D}\equiv 0$,
	\item  and  there exists the derivative $\frac{\partial v}{\partial \vec{n}_{\text{\rm in}}}\bigm|_{\partial D}\equiv 0$. 
\end{itemize}
If a non-zero function $f\in \text{\rm Hol\,}(D)$  
vanishes on a sequence  ${\tt Z}=\{{\tt z}_k\}\subset D$,
and $f$ satisfies the condition 	$ 	\int\limits_{D\setminus D_0} \log |f| \, d \nu_v<+\infty $, then 
$$
\sum\limits_{{\tt z}_k \in D\setminus D_0} v({\tt z}_k)<+\infty.
$$}
Our results allow to obtain various forms of development and generalization of results from [2]--[22], as well as new results both 
for the complex plane $\mathbb C$ and for domains in $\mathbb C_{\infty}$ such as $\mathbb D$, and another domains.
Detailed treatment of these and other results are submitted to the journal <<Математический сборник>>.

%%Пример оформления литературы:
\smallskip
\smallskip \centerline{\bf Литература}\nopagebreak

\noindent
1. {\it Ransford~Th.\/} Potential Theory in the Complex Plane --- Camb\-r\-i\-d\-ge:  Cambridge University Press,~1995.

\noindent
2. {\it  Шамоян~Ф.\,А.\/} О нулях аналитических в круге функций, рас\-т\-у\-щ\-их вблизи границы // Изв. АН Арм.ССР. Математика. 
--- 1983. --- Т.~XVIII,  №1. --- С.~15--27.

\noindent
3. {\it Ф. А. Шамоян\/} О нулях аналитических в круге функций с задан\-н\-ой мажорантой вблизи его границы // Матем. заметки. --- 2009.
---  85\,:\,2. ---  P. 300--312.

\noindent
4. {\it Шамоян~Ф.\,А., Родикова~Е.\,Г.\/} О характеризации корневых множ\-еств одного весового
класса аналитических в круге функций // Вла\-д\-и\-к\-а\-в\-к\-а\-з\-с\-к\-ий  матем. журнал. ---
2014. --- Т.~16, №~3. --- С.~64--75.

\noindent
5. {\it Е.\,Г.~Родикова\/} Факторизационное представление и описание корн\-е\-в\-ых множеств одного класса аналитических в круге функций // Сиб. электрон. матем. изв. --- 2014. --- 11. --- С.~52--63. 

\noindent
6. {\it Ф.\,А.~ Шамоян, В.\,А.~Беднаж, О.\,В.~Карбанович\/} О классах анали\-т\-и\-ч\-е\-с\-к\-их в круге функций с характеристикой Р. Неванлинны и $\alpha$-характеристикой из весовых $L^p$ пространств // Сиб. электрон. мат\-ем. изв. --- 2015. --- 12. --- С.~150--167. 

\noindent
7.  {\it Hedenmalm~H., Korenblum~B., Zhu~K.\/}  Theory of Bergman spaces. --- N.\,Y.: Springer--Verlag, Graduate Texts in Math.~2000.

\noindent
8. {\it Хабибуллин~Б.\,Н.\/} Нулевые подмножества, представление ме\-р\-о\-м\-о\-р\-ф\-н\-ых функций и характеристики Неванлинны в круге // Матем. сб. --- 2006. --- Т. 197, № 2. --- С. 117-136.  

\noindent
9.  {\it Гольдберг~А.\,А., Левин~Б.\,Я., Островский~И.\,В.\/} Целые и ме\-р\-о\-м\-о\-р\-ф\-н\-ые функции // Итоги науки и техники. Совр. про\-б\-л. матем. Фундам. напр.  --- Т.85. --- ВИНИТИ. --- М. --- 1991. --- С. 5--185.

\noindent
10. {\it Хабибуллин~Б.\,Н.\/} Полнота систем экспонент и множества единс\-т\-в\-е\-н\-н\-о\-с\-ти (издание $4^{\text{\underline{ое}}}$, дополненное). Уфа: РИЦ БашГУ. 2012

\noindent
http://www.researchgate.net/profile/Bulat{\_}Khabibullin/contributions

\noindent
11. {\it С.\,В.~Быков, Ф.\,А.~Шамоян\/} О нулях целых функций с мажора\-н\-т\-ой бесконечного порядка // Алгебра и анализ --- 2009. --- 21\,:\,6. --- С.~66--79.

\noindent
12. {\it A.~Borichev, L.~Golinskii, S.~Kupin\/} A Blaschke-type condition and its application to complex Jacobi
matrices~// Bull. London Math. Soc. --- 2009. --- 41. ---  P.~117--123.

\noindent
13. {\it S.~Favorov, L.~Golinskii\/} A Blaschke-type condition for analytic and subharmonic functions and application
to contraction operators\,// Amer. Math. Soc. Transl. --- 2009. --- 226 (2). --- P.~37–47.

\noindent
14.  {\it S.\,Yu.~Favorov, L.\,B.~Golinskii\/} Blaschke-Type Conditions for Analytic and Subharmonic Functions in
the Unit Disk~//  Local Analogs and Inverse Problems, Computational Methods and Function Theory. --- 2012. ---
 12, №1. --- P.~151--166.

\noindent
15. {\it Golinskii L., Kupin S.\/}   A Blaschke-type condition for analytic functions on finitely connected domains.
Applications to complex perturbations of a finite-band selfadjoint operator 
// J.  Math. Anal.  Appl. --- 2012. --- 389\,:\,2. --- P.~705--712.

\noindent
16. {\it S.\,Ju.~Favorov and L.\,D.~Radchenko\/} On Analytic and Subharmonic Functions in Unit Disc Growing Near a Part of the Boundary~// 
Journal of Math. Phys., Analysis, Geom. --- 2013. --- V. 9, No.~3. --- P.~304--315.

\noindent
17. {\it S.~Favorov, L.~Golinskii\/} Blaschke-type conditi\-ons in unbounded domains, generalized convexity and applications in perturbation theory~//
http://arxiv.org/pdf/1204.4283.pdf (to appear).

\noindent
18. {\it С. Ю. Фаворов, Л. Д. Радченко\/} Мера Рисса функций, субгарм\-о\-н\-и\-ч\-е\-с\-к\-их во внешности компакта~//
Математичнi Студiї. --- 2013. --- Т.~40, №~2. --- С.~149--158.

\noindent
19. {\it Хабибуллин~Б.\,Н.\/} Последовательности нулей голомофных фу\-н\-к\-ц\-ий, представление мероморфных функций и гармонические ми\-н\-о\-р\-а\-н\-ты // Матем. сб. --- 2007. --- 198\,:\,2. --- С.~121--160.

\noindent
20.  {\it Хабибуллин~Б.\,Н., Хабибуллин~Ф.\,Б., Чередникова~Л.\,Ю.\/} Под\-п\-о\-с\-л\-е\-д\-о\-в\-а\-т\-е\-л\-ь\-н\-о\-с\-ти нулей для классов голоморфных функций, их устойчивость и энтропия линейной связности. II // Алгебра и а\-н\-а\-л\-из. --- 2008. ---  20\,:\,1. --- С.~190--236.

\noindent
21.  {\it Кудашева~Е.\,Г., Хабибуллин~Б.\,Н.\/} Распределение нулей голо\-м\-о\-р\-ф\-н\-ых функций умеренного роста в единичном круге \dots
%%и представление в нём мероморфных функций 
// Матем. сб.  --- 2009. --- Т.~200, №~9. --- С.~95--126.

\noindent
22.   {\it B. N. Khabibullin\/} Generalizations of Nevanlinna’s Theorems~// Ма\-т\-е\-м\-а\-тичнi Студiї. ---Т.~34, №~2. --- 2010.

\noindent
23.  {\it Arsove~M.\, G.\/}  Functions representable as differences of subharmonic functions // Trans. Amer. Math. Soc. --- 1953. ---
 V. 75. --- P. 327--365.

\noindent
24. {\it  А. Ф. Гришин, Нгуен Ван Куинь, И. В. Поединцева\/} Теоремы о представлении $\delta$-субгармонических функций //
Biсник Харкiвського нацiонального унiверситету iменi В.\,Н.~Каразiна. Cepiя "Матема\-т\-и\-ка, прикладна математика i механiка". ---2014. ---№ 1133. ---С.~56-75; 

\noindent
{\small http://vestnik-math.univer.kharkov.ua/Vestnik-KhNU\--1133-2014-grish.pdf}

\smallskip
\noindent
{\sc Информация об авторах.} Хабибуллин Булат Нурмиевич, д.~ф.-м.~н., проф., Башкирский государственный университет (БашГУ); 
Та\-м\-и\-н\-д\-а\-р\-о\-ва (Арипова) Наргиза Рустамовна, аспирант, БашГУ
\end{document}